\documentclass{article}

\usepackage{amsmath,amssymb,amsthm}
\usepackage{xcolor}
\usepackage[colorlinks=true]{hyperref}

\title{Addendum to \textit{Pontryagin's maximum principle for dynamic systems on time scales}}

\author{Lo\"ic Bourdin\footnote{Universit\'e de Limoges, Institut de recherche XLIM, P\^ole MATHIS. UMR CNRS 7252. Limoges, France (\texttt{loic.bourdin@unilim.fr}).}
\and
Oleksandr Stanzhytskyi\footnote{Taras Shevchenko National University of Kyiv, Kyiv, Ukraine (\texttt{ostanzh@gmail.com})}
\and
Emmanuel Tr\'elat\footnote{Sorbonne Universit\'es, UPMC Univ Paris 06, CNRS UMR 7598, Laboratoire Jacques-Louis Lions, Institut Universitaire de France, F-75005, Paris, France (\texttt{emmanuel.trelat@upmc.fr}).}
}

\begin{document}

\maketitle

\begin{abstract}
This note is an addendum to \cite{BT,BKLS}, pointing out the differences between these papers and raising open questions.
\end{abstract}

\bigskip


\bigskip


\paragraph{The main differences.}
In view of establishing a time scale version of the Pontryagin Maximum Principle (PMP), the authors of \cite[Theorem 1]{BT} have developed in 2013 a strategy of proof based on the \textit{Ekeland variational principle}. This strategy was originally considered for the classical continuous case by Ivar Ekeland in his seminal paper \cite{Ekeland}.

The authors of \cite[Theorem 2.11]{BKLS} developed in 2017 a different approach, with \textit{packages of needle-like variations} and \textit{necessary conditions for an extreme in a cone}. Note that the authors of \cite{BKLS} prove moreover in \cite[Theorem 2.13]{BKLS} that the necessary conditions derived in the PMP are also sufficient in the linear-convex case.

In the sequel of this paragraph, we focus on the major pros and cons of each approach:
\begin{enumerate}
\item In \cite{BT}:
\begin{enumerate}
\item The set $\Omega$ of control constraints is assumed to be closed. This is in order to apply the Ekeland variational principle on a complete metric space.
\item There is no assumption on the time scale $\mathbb{T}$.
\end{enumerate}
\item In \cite{BKLS}:
\begin{enumerate}
\item The set $\Omega$ of control constraints is assumed to be convex, but need not to be closed.
\item The time scale $\mathbb{T}$ is assumed to satisfy \textit{density conditions} (see \cite[Definition 2.4]{BKLS}) of the kind
$$ \lim\limits_{\substack{\beta \to 0^+ \\ s+ \beta \in \mathbb{T}}} \dfrac{\mu (s+\beta)}{\beta} = 0, $$
for every right-dense points $s$, in order to guarantee that
$$
\lim\limits_{\beta \to 0^+} \; \frac{1}{\beta}  \int_{[s,s+\beta)_\mathbb{T}} x(\tau) \, \Delta \tau = x(s),
$$
for $\Delta$-integrable function $x$ and for right-dense $\Delta$-Lebesgue points $s$, even for $\beta > 0$ such that $s+\beta \notin \mathbb{T}$. Note that a discussion about this issue was provided in \cite[Section~3.1]{BT}.
\end{enumerate}
\end{enumerate}
Hence, the method developed in \cite{BKLS} allows to remove the closedness assumption done on $\Omega$ in \cite{BT}, but this is at the price of an additional assumption on the time scale $\mathbb{T}$. 

In \cite[Section 3.1]{BT}, the authors explained why other approaches (other than the Ekeland variational principle), based for instance on implicit function arguments, or on Brouwer fixed point arguments, or on separation (Hahn-Banach) arguments, fail for general time scales.

As a conclusion, a time scale version of the PMP without closedness assumption on $\Omega$ and without any assumption on the time scale $\mathbb{T}$ still remains an open challenge.

\paragraph{Additional comments on the terminal constraints.}
In \cite{BT} the authors considered constraints on the initial/final state of the kind $g (x(t_0),x(t_1)) \in \mathrm{S}$, where $\mathrm{S}$ is a nonempty closed convex set and $g$ is a general smooth function. 

In \cite{BKLS} the authors considered constraints on the initial/final state of the kind $\Phi_i (x(t_0),x(t_1)) = 0$ for $i=1,\ldots,k$, and $\Phi_i (x(t_0),x(t_1)) \leq 0$ for $i=k+1,\ldots,n$, where $\Phi_i$ are general smooth functions. 

Contrarily to what is claimed in \cite{BKLS}, the terminal constraints considered in \cite{BKLS} are only a particular case of the ones considered in \cite{BT}. Indeed, it suffices to take 
$$ g = (\Phi_1, \ldots , \Phi_k , \Phi_{k+1}, \ldots , \Phi_n )$$ 
and 
$$\mathrm{S} = \{ 0 \} \times \ldots \times \{ 0 \} \times \mathbb{R}^- \times \ldots \times \mathbb{R}^- .$$

Moreover, note that the necessary condition $-\Psi \in \mathcal{O}_\mathrm{S} ( g(x(t_0),x(t_1) ) )$ obtained in \cite[Theorem~1]{BT} encompasses both the \textit{sign condition (1)} and the \textit{complementary slackness (2)} obtained in \cite[Theorem 2.11]{BKLS}. For the sign condition, it is sufficient to recall that the orthogonal of $\mathbb{R}^-$ at a point $x \in \mathbb{R}^-$ is included in $\mathbb{R}^+$. For the complementary slackness, it is sufficient to recall that the orthogonal of $\mathrm{S}$ at $g(x(t_0),x(t_1) )$ is reduced to $\{ 0 \}$ when $g(x(t_0),x(t_1) )$ belongs to the interior of~$\mathrm{S}$.

\paragraph{Additional comments on the convexity of $\Omega$.}
The set $\Omega$ is assumed to be convex in \cite{BKLS}, while it is not in \cite{BT}. As explained in \cite[Section 3.1]{BT}, in order to apply necessary conditions of an extreme in a cone, the authors of \cite{BKLS} require that the parameters of perturbations live in intervals. As a consequence, in order to remove the convexity assumption on $\Omega$, one would need (local-directional) convexity of the set~$\Omega$ for perturbations at right-scattered points, which is 
a concept that differs from the stable $\Omega$-dense directions used in \cite{BT}. Hence, in spite of the claim done in \cite{BKLS}, the convexity assumption on $\Omega$ does not seem to be easily removable.
%
%

\paragraph{On the universal Lagrange multipliers.} 
This paragraph is devoted to providing more details on the existence of universal Lagrange multipliers claimed in \cite[page 25]{BKLS}. In the sequel, we use the notations of \cite{BKLS}, and we denote by $\mathcal{S}$ the unit sphere of $\mathbb{R}^{n+1}$.

A package $\mathrm{P}$ consists of:
\begin{itemize}
\item[-] $N \in \mathbb{N}$ and $\nu \in \mathbb{N}$;
\item[-] $\overline{\tau} = (\tau_1,\ldots,\tau_N)$ where $\tau_i$ are right-dense points of $\mathbb{T}$;
\item[-] $\overline{v} = (v_1,\ldots,v_N)$ where $v_i \in U$;
\item[-] $\overline{r} = (r_1,\ldots,r_\nu)$ where $r_i$ are right-scattered points of $\mathbb{T}$.
\item[-] $\overline{z} = (z_1,\ldots,z_\nu)$ where $z_i \in U$.
\end{itemize} 
Let $(\mathrm{P}_i)_{i \in I}$ denotes the set of all possible packages. 

Following the proof of \cite[Theorem 2.11]{BKLS}, for every $i \in I$, there exists a nonzero vector $\lambda = (\lambda_0, \ldots , \lambda_n)$ (that we renormalize in $\mathcal{S}$) of Lagrange multipliers such that:
\begin{enumerate}
\item[(i)] (1) and (2) in \cite[Theorem 2.11]{BKLS} are satisfied;
\item[(ii)] the adjoint vector $\Psi$ solution of (2.9), with the final condition (3.65) which depends on~$\lambda$, satisfies the initial condition $\Psi (t_0) = L_{x_0}$;
\item[(iii)] (4a) and (4b) in \cite[Theorem 2.11]{BKLS} are satisfied, but {only at the points contained in $\overline{\tau}$ and $\overline{r}$ respectively.}
\end{enumerate}
For every $i \in I$, the above vector $\lambda$ is not necessarily unique. Then, for every $i \in I$, we denote by $\mathrm{K}_i$ the set of all nonzero and renormalized Lagrange multiplier vectors  associated with $\mathrm{P}_i$ satisfying the above properties.

By continuity of the adjoint vector $\Psi$ with respect to the Lagrange multipliers (dependence from its final condition), we infer that $\mathrm{K}_i$ is a nonempty closed subset contained in the compact $\mathcal{S}$. This is true for every $i \in I$.

Now, let us prove that the family $(\mathrm{K}_i)_{i \in I}$ satisfies the finite intersection property. Let $J \subset I$ be a finite subset and let us prove that $\cap_{i \in J} \mathrm{K}_i \neq \emptyset$. Note that we can construct a package $\mathrm{P}$ corresponding to the union of all packages $\mathrm{P}_i$ with $i \in J$. It follows that $\mathrm{P} \in (\mathrm{P}_i)_{i \in I}$, and thus there exists a nonzero and renormalized Lagrange multiplier vector $\lambda$ associated with $\mathrm{P}$ satisfying the above properties. Since $\lambda \in \mathrm{K}_i$ for every $i \in J$, we conclude that $\cap_{i \in J} \mathrm{K}_i \neq \emptyset$.

It follows from the lemma of a centered system in a compact set that $\cap_{i \in I} \mathrm{K}_i \neq \emptyset$, and we deduce the existence of a universal Lagrange multiplier vector.

\paragraph{On the density conditions and the Cantor set.}
Contrarily to what is claimed in~\cite[Example~2.5]{BKLS}, the classical Cantor set does not satisfy the density conditions. However, generalized versions of the Cantor set (see, e.g., \cite{Turbin}) that satisfy density conditions can be constructed as follows.

Let $(\alpha_k)_{k \in \mathbb{N}}$ be a real sequence such that $0 < \alpha_k < \frac{1}{2}$ for all $k \in \mathbb{N}$, and such that $\lim_{k \to +\infty} \alpha_k = \frac{1}{2}$. Let $(A_k)_{k \in \mathbb{N}}$ be the sequence of compact subsets defined by the induction 
$$ 
A_0 = [0,1], \qquad
A_{k+1} = \mathcal{T}_k ( A_k ) \quad \forall k \in \mathbb{N},
$$
where $\mathcal{T}_k$ denotes the operator removing the open $(\alpha_k,1-\alpha_k)$-central part of all intervals. Note that the classical Cantor set corresponds to the case where $\alpha_k = \frac{1}{3}$ for every $k \in \mathbb{N}$. 

In our situation, we obtain  
$$ A_1 = [0,\alpha_0] \cup [1-\alpha_0,1], $$
$$ A_2 = \Big( [0,\alpha_1 \alpha_0] \cup [(1-\alpha_1)\alpha_0,\alpha_0] \Big) \cup \Big( [1-\alpha_0,1-(1-\alpha_1)\alpha_0)] \cup [1-\alpha_1 \alpha_0,1] \Big), $$
etc. We define the generalized Cantor set $\mathbb{T} = \cap_{k \in \mathbb{N}} A_k$. In order to prove that the time scale $\mathbb{T}$ satisfies the density conditions, from the fractal properties of $\mathbb{T}$, it suffices to prove that the density condition is satisfied at the right-dense point $0 \in \mathbb{T}$. More precisely, it is sufficient to prove that
$$ \lim\limits_{\substack{\beta \to 0^+ \\ \beta \in \mathbb{T}}} \dfrac{\mu (\beta)}{\beta} = 0. $$
Since $\mu (\beta) = 0$ for every right-dense point $\beta$, we only have to consider the case where $\beta$ is a right-scattered point of $\mathbb{T}$. In that case, one can easily see that $\frac{\mu (\beta)}{\beta} \leq \frac{1-2 \alpha_k}{\alpha_k}$ for some $k \in \mathbb{N}$ and that $k$ tends to $+\infty$ when $\beta$ tends to $0$. The conclusion follows from the fact that $\lim_{k \to +\infty} \alpha_k = \frac{1}{2}$.

%

%

\end{document}